\documentclass[11pt]{amsart}
\allowdisplaybreaks[4]
\usepackage{amssymb,color}
\usepackage[square,compress,comma, numbers,sort]{natbib}

\usepackage[colorlinks=true, citecolor=red, linkcolor=blue]{hyperref}
\usepackage{amsthm}
\usepackage{mathtools}

\DeclareMathOperator{\argmax}{argmax}
\DeclareMathOperator{\var}{Var}
\DeclareMathOperator{\corr}{Corr}

\definecolor{c20}{rgb}{0.,0.7,0.}
\definecolor{c30}{rgb}{0.,0.,1.}
\definecolor{c40}{rgb}{1,0.1,0.7}
\definecolor{c50}{rgb}{1,0,0}
\definecolor{c60}{rgb}{0,0.9,0.1}

\newcommand{\ve}{\varepsilon}

\newcommand{\E}[1]{\mathbb{E} \left(#1\right)}

\newcommand{\pk}[1]{\mathbb{P} \left(#1 \right) }

\newcommand{\R}{\mathbb{R}}

\newcommand{\ldot}{,\ldots,}

\newcommand{\limit}[1]{\lim_{#1 \to   \infty}}

\newcommand{\BQN}{\begin{eqnarray}}
\newcommand{\EQN}{\end{eqnarray}}
\newcommand{\BQNY}{\begin{eqnarray*}}
\newcommand{\EQNY}{\end{eqnarray*}}

\def\nncol#1{\textcolor{black}{#1}}

\def\K1#1{\textcolor{cyan}{#1}}
\def\K1#1{#1}

\def\bqny#1{ \nncol{ \begin{eqnarray*} #1 \end{eqnarray*}}}
\def\bqn#1{ \nncol{ \begin{eqnarray} #1 \end{eqnarray}}}

\newcommand{\BS}{\begin{sat}}
\newcommand{\ES}{\end{sat}}
\newcommand{\BT}{\begin{theo}}
\newcommand{\ET}{\end{theo}}
\newcommand{\BK}{\begin{korr}}
\newcommand{\EK}{\end{korr}}

\newcommand{\BD}{\begin{de}}
\newcommand{\ED}{\end{de}}

\newcommand{\BIT}{\begin{itemize}}
\newcommand{\EIT}{\end{itemize}}
\newcommand{\BDI}{\begin{description}}
\newcommand{\EDI}{\end{description}}

\newcommand{\BRM}{\begin{remark}}
\newcommand{\ERM}{\end{remark}}

\newcommand{\BEL}{\begin{lem}}
\newcommand{\EEL}{\end{lem}}

\newtheorem{theo}{Theorem}[section]
\newtheorem{sat}[theo]{Proposition}
\newtheorem{de}[theo]{Definition}
\newtheorem{lem}[theo]{Lemma}

\newtheorem{korr}[theo]{Corollary}
\newtheorem{remark}[theo]{Remark}

\newcommand{\prooftheo}[1]{ \textbf{Proof of Theorem} \ref{#1} }

\newcommand{\COM}[1]{}

\newcommand{\QED}{\hfill $\Box$ \\}

\def\rw{\rightarrow}

\def\IF{\infty}

\def\eHH#1{\textcolor{c50}{#1}}

\def\eHH#1{#1}

\def\rred#1{\textcolor{black}{#1}}

\def\vf{\sigma^2}

\newcommand{\mQ}{\mathbb{C}}

\begin{document}

\title
{Extremes of reflecting Gaussian processes on discrete Grid}

\author{Krzysztof D\c{e}bicki}
\address{Krzysztof D\c{e}bicki, Mathematical Institute, University of Wroc\l aw, pl. Grunwaldzki 2/4, 50-384 Wroc\l aw, Poland}
\email{Krzysztof.Debicki@math.uni.wroc.pl}

\COM{
\author{Enkelejd  Hashorva}
\address{Enkelejd Hashorva, Department of Actuarial Science 
	University of Lausanne,\\
	UNIL-Dorigny, 1015 Lausanne, Switzerland}
\email{Enkelejd.Hashorva@unil.ch}
}

\author{{Grigori Jasnovidov}}
\address{Grigori Jasnovidov, St. Petersburg Department
of Steklov Mathematical Institute
of Russian Academy of Sciences, St. Petersburg, Russia;
University of Lausanne, Lausanne, Switzerland
}
\email{griga1995@yandex.ru}

\bigskip

 \maketitle
\bigskip
{\bf Abstract:} For $\{X(t), t \in G_\delta\}$ a centered Gaussian process
with stationary increments 
on a discrete grid
$G_\delta=\{0,\delta,2\delta, ...\}$, where $\delta>0$,
we investigate
the stationary reflected process
\[Q_{\delta,X}(t) = \sup\limits_{s\in [t,\IF)\cap G_\delta}\big(
X(s)-X(t)-c(s-t)\big), \ t \in G_\delta
\]
with $c>0$.
We
derive the exact asymptotics of
$\pk{\sup\limits_{t\in [0,T]\cap G_\delta} Q_{\delta,X}(t)>u}$
and
$\pk{\inf\limits_{t\in [0,T]\cap G_\delta} Q_{\delta,X}(t)>u}, $
as $u\to\infty$, with $T>0$.
It appears that $\varphi=\lim_{u\to\infty} \frac{\sigma^2(u)}{u}$ determines the asymptotics, leading to three
qualitatively different scenarios:
$\varphi=0$, $\varphi\in(0,\infty)$ and $\varphi=\infty$.
\\



\section{Introduction}\label{s.Intro}
For  $X(t), t\ge 0 $ a centered Gaussian process  with a.s. continuous sample paths,
stationary increments and
variance function $\sigma^2(t):=\var\left(X(t)\right)$ such that $\sigma^2(0)=0$,
consider the {\it reflected} (at $0$) process
\begin{equation}
\label{prob.0}
\hat{Q}_{X}(t)=X(t)-ct+\max\left( \hat{Q}_{X}(0),-\inf_{s\in[0,t]}(X(s)-cs)   \right),\ t\ge0,
\end{equation}
where $c>0$.

Due to its relation with the solution of
the {\it Skorokhod problem} that describes the dynamics of
the {\it buffer content process}
in a {\it fluid} queue fed by $X$ and emptied at rate $c$,
properties of $\hat{Q}_{X}$ have been investigated for a wide class of processes $X$,
see e.g. \cite{Nor94,HuP99,Piterbarg2001,DeK} and references therein.
Within this framework, the process $\hat{Q}_{X}(t)$
is called in the literature as the {\it storage process} \cite{Piterbarg2001, AlS04}.

Distributional properties of the unique stationary solution of \eqref{prob.0}, which has the following representation
\begin{equation}
\label{eq:stat1}
Q_{X}(t)=\sup_{t\le s}\left(X(s)-X(t)-c(s-t)\right),
\end{equation}
were intensively analyzed.
In particular, the exact asymptotics of $ \pk{Q_X(0)>u}$,
as $u\to\infty$,
was derived for a wide class of Gaussian processes $X$;
see, e.g., \cite{HuP99, HuP04,DeK, DI2005} and references therein.
An important direction of research on extremal behaviour of
the process $Q_X(t),t\ge 0$ was initiated by Piterbarg \cite{Piterbarg2001},
where asymptotics of $\sup_{t\in[0,T]} Q_{B_\alpha}(t)$
was derived, with $B_\alpha(t)$ a fractional Brownian motion with Hurst parameter $\alpha\in(0,1)$.
Notably, the following asymptotic equivalence
\[
\pk{Q_{B_\alpha}(0)>u}\sim \pk{\sup_{t\in[0,T]} Q_{B_\alpha}(t)>u}
\]
for $\alpha >1/2$, called in \cite{AlS04} the {\it Piterbarg property}, took particular interest
and was further observed for more broad class of processes $X$ \cite{AlS04, KrzysPeng2015, DeK}.
Complementary, in \cite{DeK} it was shown that if $\alpha>1/2$, then the process $Q_{B_\alpha}$ also possesses
{\it strong Piterbarg} property, i.e. for all $T>0$
\[
\pk{Q_{B_\alpha}(0)>u}\sim \pk{\inf_{t\in[0,T]} Q_{B_\alpha}(t)>u}
\]
as $u\to\infty$.

From the point of view of the stochastic modelling or simulation techniques,
discrete-time models
frequently appear to be more natural.
However, despite of its relevance in modelling of, e.g., queueing systems, much less is known
on distributional properties of the discrete counterpart of (\ref{eq:stat1}), i.e.,
\bqn{\label{storage_ruin}
Q_{\delta,X}(t) = \sup\limits_{s\in [t,\IF)\cap G_\delta}\big(
X(s)-X(t)-c(s-t)\big), \ t \in G_\delta,
}
where
$G_\delta= \{0,\delta, 2 \delta \ldot \}$.
A notable exceptions are recent works 
\cite{Piterbargdiscrete,secondproj},
where the exact asymptotics of $\pk{Q_{\delta,B_\alpha}(0)>u}$, as $u\to\infty$, was derived.

In this contribution we extend the findings of
\cite{Piterbargdiscrete,HashorvaGrigori,secondproj} to a more general class of
Gaussian processes with stationary increments and derive the
exact asymptotics of
\bqn{\label{psi_def_0}
\psi_{T,\delta}^{\sup}(u) :=
\pk{\sup\limits_{t\in [0,T]\cap G_\delta} Q_{\delta,X}(t) >u}, \quad
\psi^{\inf}_{T,\delta}(u) :=
\pk{\inf\limits_{t\in [0,T]\cap G_\delta}
Q_{\delta,X}(t)>u},
}
as $u\to\infty$, for $T> 0$ and $\delta>0$, complementing results for
continuous time given in \cite{Piterbarg2001,DeK}.

It appears that the influence of the grid size $\delta$
on the asymptotics of (\ref{psi_def_0}) strongly depends on
the value of
\bqn{\label{varphi}
\varphi:=\lim_{u\rw\IF}\frac{\sigma^2(u)}{u}\in[0,\infty],
}
leading to three scenarios: $\varphi=0$,  $\varphi\in(0,\infty)$ and  $\varphi=\infty$.
The case $\varphi=\infty$ leads to the same asymptotics as its continuous-time counterpart,
which reflects the {\it long-range dependance property} of $X$ when its variance $\vf$ is
superlinear. On the other hand, the asymptotics for the case $\varphi=0$
strongly depends on the grid size and is asymptotically negligible with respect to its continuous
counterpart. The third case $\varphi\in(0,\infty)$ needs particularly precise analysis and
leads to the asymptotics which differs with the continuous-time case only by a constant (depending on $\delta$).

{\it Outline of the paper.}  We introduce notation and assumptions on the process $X$ in Section \ref{s.preliminary}.
Then, in Theorem \ref{TH1} we derive exact asymptotics of
$\pk{Q_{\delta,X}(0)>u}$ as $u\to\infty$. Main results of this contribution are presented in Section \ref{s.Main}.
In Section \ref{s.Examples} we illustrate the main findings of this paper by analysis of two important classes
of Gaussian processes, i.e. fractional Brownian motions and Gaussian integrated processes.
Section \ref{s.proofs} consists of the proofs of the results derived in this paper.

\section{Notation and preliminary results}\label{s.preliminary}

Let $X(t), t\in \R$ be a centered Gaussian process with stationary increments, as introduced in Section \ref{s.Intro}.
Suppose that\\
\\
{\bf A}:  $\vf$ is regularly varying at $\IF$ with
index $2\alpha\in(0,\eHH{2)}$
and $\vf(t)$ is twice continuously differentiable
for any $ t\in (0,\IF)$.
Further, the first and second derivatives of $\vf$
are ultimately monotone.\\
\\
We note that it follows straightforwardly from {\bf A} that
$\varphi=0$ if $\alpha<1/2$ while for $\alpha>1/2$ we have $\varphi=\infty$.
It appears that case $\alpha=1/2$, i.e. when $\vf$ is asymptotically close to a linear function,
needs particularly precise analysis, for which the following condition 
is a tractable assumption:\\
\\
{\bf B}: If $\vf$ satisfies {\bf A} with $\alpha = \frac{1}{2}$, then
$\varphi>0.$
\\
\\
Condition {\bf B} excludes the cases when $\alpha=1/2$ but $\varphi=0$.
For this scenario we were able to give only partial results.
We refer to Remark \ref{remark} for the discussion of
the extension of the results derived in Section \ref{s.Main}
to the case $\alpha=1/2, \varphi=0$ under an additional constrain on the variance
function $\sigma^2$ of $X$.

Conditions {\bf A} and {\bf B} are satisfied for a wide class of Gaussian processes with stationary increments,
including family of fractional Brownian
motions and integrated stationary Gaussian processes; see Section \ref{s.Examples} for details.
We note that quantity $\varphi$ already appeared in
\cite{DI2005}, where it was observed that the form of the asymptotic behavior, as $u\to\infty$,
of $\pk{Q_X(0)>u}$
introduced in (\ref{eq:stat1})
is
determined by the value of $\varphi$.
Let
\bqn{\label{H_inf} \ \ \ \ \ \
\mathcal H_\xi(M) = \E{\sup\limits_{t\in M}e^{\sqrt 2 \xi(t)-\var(\xi(t))}}
\in (0,\IF),
\quad
\mathcal G_\xi(M) = \E{\inf\limits_{t\in M}
e^{\sqrt 2 \xi(t)-\var(\xi(t))}}\in (0,\IF),
}
where $M$ is a compact subset of $\R$
and $\xi(t), \ t \in \R$ 
is a Gaussian process
with stationary increments and a.s. continuous
sample paths.
Then, we define {\it Pickands constant} by
$$  \mathcal{H}_{\xi}^{\delta} = \limit{S} \frac{
 \mathcal H_\xi([0,S]\cap\delta\mathbb Z)}{S}, \ \ \ \ \ \ \ \delta\ge 0,
$$
where we set
$\delta \mathbb{Z}=\R^+$ if $\delta=0$. We refer to \cite{DE2002} for properties of
$\mathcal{H}_{\xi}^{0}$. 
In Lemma
\ref{lemma_finite_discrete_pick_const}
we prove that for $\delta>0$
 it is sufficient to suppose that
$\xi$ satisfies {\bf A} to claim that $\mathcal H_\xi^\delta
\in (0,\IF)$.
Later on for $\delta = 0$ we simply write $\mathcal{H}_{\xi}$
instead of $\mathcal{H}_{\xi}^{0}$.

Next, let us recall the findings of
\cite{DI2005}[Proposition 2]  (see also \cite{KrzysPeng2015} [Theorems 3.1-3.3])
for the asymptotics of
$\pk{Q_X(0)>u}$, as $u\to\infty$, which will be a useful benchmark for the results
derived in the next section.
Let  $\overleftarrow{\sigma}(t), t\geq 0$
stand for the asymptotic inverse function of $\sigma$,
i.e., $\overleftarrow\sigma(x) = \inf \{y \in [0,\IF):
f(y)>x\} $
 (for details and properties of the asymptotic inverse
 functions see, e.g.,  \cite{LeadbetterExtremes})
 and  let
\begin{eqnarray}\label{theta}
t_* = \frac{\alpha}{c(1-\alpha)}, \ \ \ \
m(u)=\inf\limits_{t>0}\frac{u(1+ct)}{\sigma(ut)}, \ \ \ \
\Delta(u)=
\begin{cases}
\overleftarrow{\sigma}
\left(\frac{\sqrt{2}\sigma^2(ut_*)}{u(1+ct_*)}\right),
& \varphi \notin(0,\IF) \\
1,&  \varphi \in (0,\infty).
\end{cases}
\end{eqnarray}
Let for $X$ such that $\varphi \in (0,\IF)$,
\bqn{\label{eta}
\eta(t) =
\frac{c\sqrt 2}{\varphi}X(t), \ \ \ \ \  t\ge 0.
}
As shown in \cite{DI2005}[Proposition 2], if 
$\vf$ is regularly varying at $0$ with index $2\alpha_0\in(0,2]$ and
\textbf{A} is satisfied, then
for $\hat{Q}_X$ defined in (\ref{prob.0}), we have
\bqn{\label{cont_simpletheo}
\pk{ \hat{Q}_X(0)> u }\sim f(u)\Psi(m(u))\times
\begin{cases}
\mathcal{H}_{B_\alpha}, & \varphi = \IF\\
\mathcal{H}_{\eta}, & \varphi \in (0,\IF)\\
\mathcal{H}_{B_{\alpha_0}}, & \varphi = 0
\end{cases}
, \quad u\to \IF
,}
where $\Psi$ is the survival function of a standard Gaussian random
variable and
\rred{\BQN\label{AB}
f(u) = \sqrt{\frac{2\pi A}{B}}
\frac{u}{m(u)\Delta(u)}, \ \ \ \
	A=\frac{1}{(1-\alpha)t_*^{\alpha}}, \quad
	B=\frac{\alpha}{t_*^{\alpha+2}}.
	\EQN}

The following result establishes the asymptotics of $\pk{ Q_{\delta,X}(0)> u}$ for
 $\delta>0$ as $u \to \IF$.
It generalizes the findings of \cite{Piterbargdiscrete,secondproj},
where the special case of $X$ being a fractional Brownian motion was analyzed.

\BT\label{TH1}  Let $X(t), t\ge 0 $ be a centered Gaussian process
with continuous trajectories and stationary increments satisfying
{\bf A, B}.
Then, for $\delta>0$, as $u \to \IF$, it holds that
\bqn{\label{TH1_claim_i}
\pk{ Q_{\delta,X}(0)> u} \sim \Psi(m(u))\times
\begin{cases}
\frac{\sqrt{2\pi\alpha}u}{\delta c(1-\alpha)^{3/2} m(u)},
&\varphi =0 	\\
\mathcal{H}_{\eta}^{\delta}f(u),
&\varphi \in (0,\IF) 	\\
\mathcal{H}_{B_\alpha}f(u),
&\varphi =\IF.
\end{cases}
}	
\ET

\BRM Providing that both  the asymptotics \eqref{cont_simpletheo} and
 \eqref{TH1_claim_i} hold, we have that
\bqny{
\lim_{u\to\infty}\pk{ Q_{\delta,X}(0)> u|Q_{X}(0)>u} =
\begin{cases}
0,&  \varphi = 0,\\
 \frac{\mathcal{H}_{\eta}^{\delta}} {\mathcal{H}_{\eta}},& \varphi \in (0,\IF),\\
 1,& \varphi =\IF.
\end{cases}
}
\ERM

\section{Main Results}\label{s.Main}
In this section we derive the exact asymptotics of
\bqn{\label{psi_def}
\psi_{T,\delta}^{\sup}(u) :=
\pk{\sup\limits_{t\in [0,T]_\delta} Q_{\delta,X}(t) >u}
 \quad \text{and} \quad
\psi^{\inf}_{T,\delta}(u) :=
\pk{\inf\limits_{t\in [0,T]_\delta}
Q_{\delta,X}(t)>u}
,}
as $u\to\infty$,
for $T>0$ and $\delta>0$,
where for any real $a<b$ and positive $\delta$
$$[a,b]_\delta = [a,b]\cap \delta \mathbb{Z}.$$


We begin with the asymptotics of $\psi_{T,\delta}^{\sup}(u)$ as $u\to\infty$.
Let in the following $[\cdot]$ stand for the integer part.
\BT\label{TH2}
Let $X(t), t\ge 0 $ be a centered Gaussian process
with continuous trajectories and stationary increments satisfying
{\bf A-B}.
Then for $\delta>0$ as $u \to \IF$ it holds that
\bqny{
\psi_{T,\delta}^{\sup}(u) \sim
\Psi(m(u)) \times
\begin{cases}
(1+[\frac{T}{\delta}])\frac{\sqrt{2\pi\alpha}u}{\delta
c(1-\alpha)^{3/2} m(u)},
&\varphi =0 	\\
\mathcal{H}_\eta([0,T]_\delta)\mathcal{H}_{\eta}^{\delta}f(u),
&\varphi \in (0,\IF) 	\\
\mathcal{H}_{B_\alpha}f(u),
&\varphi =\IF.
\end{cases}
}	
\ET

\BRM \label{remark}
It follows straightforwardly from the proof of Theorem \ref{TH1} and \ref{TH2}
that condition {\bf B} can be relaxed a bit.
Namely,  if
 $\varphi =0$ and  $\alpha = \frac{1}{2}$ and
for $\kappa = \sqrt{c\!\!\inf\limits_{t\in\{\delta,2\delta,...\}}\sigma(t)}
-\ve$, with sufficiently small $\ve>0$,
\bqn{\label{theo1_case_i_assumption_sigma}
\sigma(u) \le \kappa \frac{\sqrt u}{\ln^{1/4} u},  \ \ \  \
u \to \IF}
then both Theorem \ref{TH1} and \ref{TH2} hold.
If  $\varphi=0$ and $\alpha = 1/2$ in the Theorem \ref{TH2} and
\eqref{theo1_case_i_assumption_sigma} does not hold, then
it follows from its proof that
\eqref{TH1_claim_i} reduces to the upper bound.
\ERM

Next we analyze the asymptotical behaviour of
$\psi^{\inf}_{T,\delta}(u)$, as $u\to\infty$.
\BT\label{THinf}
Let $X(t), t\ge 0 $ be a centered Gaussian process with continuous
trajectories and stationary increments satisfying
{\bf A-B}. Then for $\delta>0$ as $u \to \IF$ it holds that
\bqny{
\psi_{T,\delta}^{\inf}(u) \sim f(u)\Psi(m(u))\times
\begin{cases}
\mathcal G_\eta([0,T]_\delta)\mathcal{H}_{\eta}^{\delta},&
\varphi \in (0,\IF)\\
\mathcal{H}_{B_\alpha},&
\varphi = \IF.\end{cases}
}	
\ET

If $\varphi = 0$, then for the non-degenerated scenario (when set
$[0,T]_\delta$ consists of more than 1 element, i.e, for $T\ge \delta$)
it seems  difficult to derive even
logarithmic asymptotics of
$\psi_{T,\delta}^{\inf}(u)$.
One can argue that $\psi_{T,\delta}^{\inf}(u)$
is exponentially smaller than
$\pk{ Q_{\delta,X}(0)> u}$ in this case, as $u\to\infty$.
We have the following proposition giving an
upper bound for $\psi_{T,\delta}^{\inf}(u)$.

\begin{sat} \label{proposition1} If $\varphi = 0$ and for some
 $\ve>0$
\bqn{\label{additional_condition_proposition}
\sigma(u) \le \frac{\sqrt u}{\ln^{1/4+\ve}u}, \ \ \ \ \ \ \ u \to \IF
,}
then for $T\ge \delta$ with any
$\widetilde \mQ <\frac{1+ct_*}{2t_*^{2\alpha}}
 \sup\limits_{t\in [0,T]_\delta}\sigma(t)$ it holds that
\bqny{
\psi_{T,\delta}^{\inf}(u) \le \Psi(m(u))\Psi\left(\widetilde\mQ
\frac{u}{\sigma^2(u)}\right),  \ \ \ \ u \to \IF
.}
\end{sat}
\BRM
Comparing the asymptotics in Theorems \ref{TH1}, \ref{TH2} and \ref{THinf}
for case $\alpha>1/2$ we observe that the process $Q_{\delta,X}$ possesses the so-called strong Piterbarg property,
that is
\bqny{
\psi_{T,\delta}^{\sup}(u)
\sim \pk{Q_{\delta,X}(0)>u} \sim
\psi_{T,\delta}^{\inf}(u), \ \ \ \ \ \ \
T\ge 0, \ u \to \IF.
}
The analogous property was observed for the continuous-times
analog of $Q_{\delta,X}$; see \cite{DeK}.
\ERM

\section{Examples}\label{s.Examples}
We illustrate the findings of this contribution by application of Theorems \ref{TH2} and
\ref{THinf} to the family of fractional Brownian motions and Gaussian integrated processes.

\subsection{Fractional Brownian motion.}
Let $$C_H = \frac{c^H}{H^H(1-H)^{1-H}}, \ \ D_H =
 \frac{\sqrt{2\pi} H^{H+1/2}}{c^{H+1}(1-H)^{H+1/2}},
 \ \ E_H = \frac{2^{\frac{1}{2}-\frac{1}{2H}}\sqrt \pi}
{H^{1/2}(1-H)^{1/2}}.$$
Applying Theorems \ref{TH1}, \ref{TH2} and \ref{THinf}
 for $X(t)=B_H(t)$ being a standard fractional Brownian motion with Hurst parameter $H\in (0,1)$, we obtain the following results.
\begin{korr}\label{korr1}
 As $u \to \IF$ it holds that
\bqny{
\pk{Q_{\delta,B_H}(0)
>u} \sim
\begin{cases}
\frac{D_Hu^H}{\delta }
\Psi(C_Hu^{1-H}), & H<1/2
\\
\mathcal{H}^{2c^2\delta}_{B_{1/2}}  e^{-2cu}, & H=1/2\\
\mathcal{H}_{B_H}E_H
(C_Hu^{1-H})^{1/H-1}\Psi(C_Hu^{1-H}), & H>1/2.
\end{cases}
} \end{korr}

\begin{korr}\label{korr2} For $T,\delta>0$ as $u \to \IF$ it holds that
\bqny{
\pk{\sup\limits_{t\in [0,T]_\delta}Q_{\delta,B_H}(t)> u } \sim
\begin{cases}
(1+[\frac{T}{\delta}])\frac{D_H u^H}{\delta }
\Psi(C_Hu^{1-H}), & H<1/2\\
\mathcal{H}_{B_{1/2}}([0,2c^2T]_{2c^2\delta})
\mathcal{H}^{2c^2\delta}_{B_{1/2}}  e^{-2cu}, & H=1/2\\
\mathcal{H}_{B_H}E_H
(C_Hu^{1-H})^{1/H-1}\Psi(C_Hu^{1-H}), & H>1/2.
\end{cases}
}
and
%
%
%
%
\bqny{
\pk{\inf\limits_{t\in [0,T]_\delta}Q_{\delta,B_H}(t)> u } \sim
\begin{cases}
\mathcal G_{B_{1/2}}([0,2c^2T]_{2c^2\delta})
\mathcal{H}_{B_{1/2}}^{2c^2\delta}
e^{-2cu},& H=1/2\\
\mathcal{H}_{B_H}E_H
(C_Hu^{1-H})^{1/H-1}\Psi(C_Hu^{1-H}),&
H>1/2.\end{cases}
}\end{korr}

Note that Corollary \ref{korr1} intersects with the results
 in \cite{secondproj,Piterbargdiscrete,HashorvaGrigori}
while Corollary \ref{korr2} provides a discrete counterpart
of Theorems 5-7 in \cite{Piterbarg2001} and Theorem 1 in
\cite{DeK}, respectively.

\subsection{Gaussian integrated processes.}
For a stationary centered
Gaussian process with a.s. continuous sample paths $\zeta(s), \ s\ge 0$
define the integrated process by
\bqn{\label{Z_integrated}
Z(t) = \int\limits_0^t\zeta(s)ds, \ \ \ \ \ t \ge 0.
}
This process is also Gaussian, has a.s. continuous
 sample paths and stationary
increments. In what follows
we consider two classes
of processes $Z$,
which differ by property of
the correlation
function
$R(t) := \E{\zeta(0)\zeta(t)}$
 of $\zeta$ as $t\to\infty$.
\\ \\
\emph{SRD case.}
Following, e.g.,  \cite{Debicki2002} (see also
\cite{Debicki1995}), we impose the
following conditions on the correlation of $\zeta$:
\\
\textbf{S1}: $R(t) \in C([0,\IF)), \ \lim\limits_{t \to \IF}tR(t) = 0$;\\
\textbf{S2}: $\int\limits_0^tR(s)ds>0$ for all $t\in (0,\IF]$;\\
\textbf{S3}: $\int\limits_0^\IF t^2|R(t)|dt<\IF$.\\
The above assertions imply the existence of the
first and second derivatives of $\sigma^2_Z(t)=Var (Z(t))$ 
and establish the asymptotic behavior
of $\sigma^2_Z(t)$ at $\IF$ (see e.g., Remark 6.1 in \cite{Debicki2002}):
$$\sigma^2_Z(t) = \frac{2}{G}t-2D+o(t^{-1}), \ \ \ \ \ \ \  t \to \IF,$$
where $G = 1/\int\limits_{0}^\IF R(t)dt$ and $D =\int\limits_{0}^\IF tR(t)dt$.
Thus, $\sigma^2_Z$ satisfies \textbf{A-B} with $\alpha = 1/2$ and
$\varphi=\frac{2}{G}>0$. Hence
applying Theorems \ref{TH1}, \ref{TH2} and \ref{THinf}, the following corollary holds.
\begin{korr}\label{korr_integrated_simple}
Suppose that $\zeta$ satisfies \textbf{S1-S3}. Then for $T\ge 0$
and $\delta>0$
as $u \to \IF$
\bqny{
\pk{Q_{\delta,Z}(0)
>u}  &\sim&
\mathcal A \mathcal H^\delta_{\xi}   e^{-cGu} ,  \\
\pk{\sup\limits_{ t \in [0,T]_\delta}
Q_{\delta,Z}(t)>u }
&\sim& \mathcal A \mathcal H _{\xi} ([0,T]_\delta)
 \mathcal H^\delta_{\xi}   e^{-cGu},\\
\pk{\inf\limits_{ t \in [0,T]_\delta}
Q_{\delta,Z}(t)>u }
&\sim& \mathcal A \mathcal G _{\xi} ([0,T]_\delta)
 \mathcal H^\delta_{\xi}   e^{-cGu},
}
where $\mathcal A = \frac{1}{c^2Ge^{c^2G^2D}}$ and
$\xi(t) = cGZ(t)/\sqrt 2$.
\end{korr}

Note that the first asymptotics in Corollary \ref{korr_integrated_simple}
 differs from its
continuous-time analog (Theorem 5.1 in \cite{Debicki2002})
only by the corresponding Pickands constants.
\\ \\
\emph{LRD case.} Following, e.g., \cite{HuP04,DI2005} we characterize LRD case by
the following assumptions on the covariance function $R(t)$ of the process $\zeta$:\\
\textbf{L1}: $R(t)$ is a continuous strictly positive function for $t\ge 0$;
\\
\textbf{L2}: $R(t)$
is regularly varying at $\IF$ with index $2\alpha - 2$,
$\alpha\in (1/2,1)$.

Under the above assumptions,
by Karamata's theorem, $\sigma_Z^2$ is regularly
varying at $\IF$ with index $2\alpha$.
Since $2\alpha>1$ we are in
$\varphi = \IF$ scenario.
Hence,
applying Theorems \ref{TH1}, \ref{TH2} and
\ref{THinf} we  obtain the following result,
which shows that $Q_{\delta,Z}(t)$ possesses the strong
Piterbarg property.
\begin{korr}
Suppose that $\zeta$ satisfies
\textbf{L1-L2}. Then
as $u \to \IF$ it holds that
\bqny{
\pk{\inf \limits_{t\in [0,T]_\delta} Q_{\delta,Z}(t) > u }
\sim
\pk{
Q_{\delta,Z}(0)>u}
&\sim& \pk{\sup \limits_{t\in [0,T]_\delta} Q_{\delta,Z}(t) > u }
\\&\sim&    \mathcal{H}_{B_{\alpha}} f(u)\Psi(m(u)),
}
where $m(u)$ and $f(u)$ are defined in \eqref{theta} and \eqref{AB},
respectively.
\end{korr}

\section{Proofs}\label{s.proofs}
In this section we give proofs of all the results presented in this contribution.
Hereafter, denote by $\mathbb{C}$, $\mathbb{C}_i, \ i=1,2,3,\dots$ positive constants that may differ from line to line and $\overline{X}:=\frac{X}{\sqrt{Var(X)}}$ for any nontrivial random variable $X$.
\COM{
Moreover, $f(u)\sim g(u), u\rw\IF$ means that $\lim_{u\rw\IF}\frac{f(u)}{g(u)}=1.$  In our proofs, multiple limits appear. We shall write for instance
$$b(u,S,\epsilon)\sim a(u), \quad u\rw\IF, S\rw\IF, \epsilon\rw 0 $$
to mean that
$$\lim_{\epsilon\rw 0}\lim_{S\rw\IF}\lim_{u\rw\IF}\frac{b(u,S,\epsilon)}{a(u)}=1.$$
}
For any $u>0$ we have
$$\pk{Q_{\delta,X}(0)>u}=\pk{\sup_{t \in G_\delta }(X(t)-ct)>u}=\pk{\sup_{t\in G_{\delta/u }}X_u(t)>m(u)},$$
where $m(u)$ is defined in \eqref{theta} and
$$X_u(t)=\frac{X(ut)}{u(1+ct)}m(u).$$
Denote by $\sigma_{X_u}^2$ the variance function of $X_u(t), \ t\ge 0$.
In the next lemma we focus on asymptotic properties of
the variance and correlation functions of $X_u(t)$; we refer to, e.g., \cite{KrzysPeng2015} for the proof.

\BEL\label{L1}
Suppose that {\bf A} is satisfied.
 For $u$ large enough the maximizer $t_u$ of $\sigma_{X_u}$ is unique and
$t_u\rw t_*=\frac{\alpha}{c(1-\alpha)}$ as $u\to \infty$. Moreover,  for $\delta_u>0$
satisfying $\lim\limits_{u\rw\IF}\delta_u=0$ ($A$, $B$
are defined in \eqref{AB})
\BQNY\label{V4}
\lim_{u\rw\IF}\sup_{t\in (t_u-\delta_u, t_u+\delta_u)\setminus \{t_u\}}\left|\frac{1-\sigma_{X_u}(t)}{\frac{B}{2A}(t-t_u)^2}-1\right|=0
\EQNY
and (recall, $\sigma^2$ is the variance of $X$)
 \BQNY\label{cor0}
\lim_{u\rw \IF}\sup_{s\neq t, s,t \in (t_u-\delta_u, t_u+\delta_u)}\left|\frac{1-Cor\left(X(us), X(ut)\right)}{\frac{\sigma^2(u|s-t|)}{2\sigma^2(ut_*)}}-1\right|=0.
\EQNY
\EEL

By Lemma \ref{L1}
we have that $t_u$ is the unique minimizer of $
\frac{u(1+ct)}{\sigma(ut)}$ for large $u$ and hence by Potter's theorem
(Theorem 1.5.6 in \cite{BI1989}) we obtain useful in the following proofs
asymptotics of $m(u)$ 
\bqn{\label{asympt_m(u)}
m(u) = \frac{u(1+ct_u)}{\sigma(ut_u)} \sim \frac{u(1+ct_*)}{\sigma(ut_*)}
\cdot\frac{\sigma(ut_*)}{\sigma(ut_u)}
\sim
\frac{u(1+ct_*)}{t_*^{\alpha}\sigma(u)},\quad u \to \IF
.}

Observe that
\bqny{
\psi_{T,\delta}^{\sup}(u) =
\pk{\sup\limits_{t\in [0,T/u]_{\delta/u}
,t\le s\in G_{\delta/u}}Z_u(t,s)>m(u)},
}
where
\bqn{\label{Z_u_def}
Z_u(t,s) = \frac{X(us)-X(ut)}{u(1+c(s-t))}m(u).
}
Notice that for the variance $\sigma_{Z_u}^2$ of $Z_u$ it holds, that
$\sigma_{Z_u}^2(s,t) = \sigma_{X_u}^2(s-t)$
and for correlation $r_{Z_u}$ we have
for $\delta_u>0$ satisfying $\lim_{u\rw\IF}\delta_u=0$
(Lemma 5.4. in \cite{KrzysPeng2015})
\bqn{\label{r_Z_u_formula}
\limit{u}\sup\limits_{|t-t_1|<\delta_u,s-t,s_1-t_1\in (-\delta_u+t_u,t_u+\delta_u)
,(s,t)\neq(s_1,t_1)}
\Big| \frac{1-r_{Z_u}(s,t,s_1,t_1)}{\frac{\sigma^2(u|s-s_1|)+\sigma^2(u|t-t_1|)}
{2\sigma^2(ut_*)}}-1\Big| = 0.
}
To the rest of the paper we suppose that
\bqny{
\delta_u =
\begin{cases}
u^{-1/2}\ln u, & \varphi<\IF\\
 u^{-1}\ln(u)\sigma(u), &\varphi = \IF
\end{cases}
}
and set
\bqn{\label{I_t_u_def}
I(t_u) =
 G_{\frac{\delta}{u}}\bigcap (-\delta_u+t_u,t_u+\delta_u)
}
for $u>0$.
The following lemma allow us to extract the main area contributing in the asymtotics of
$\psi_{T,\delta}^{\sup}(u), \psi_{T,\delta}^{\inf}(u)$ and
$\pk{Q_{\delta,X}(0)>u}$ as $u \to \IF$.
\begin{lem}\label{the_main_area} For any $T\ge 0$ it holds that,
as $u \to \IF$,
\bqny{
\psi_{T,\delta}^{\sup}(u) &\sim& \pk{\sup\limits_{t\in [0,T/u]_{\delta/u}
,s\in I(t_u)}Z_u(t,s)>m(u)}\\
\psi_{T,\delta}^{\inf}(u) &\sim& \pk{\inf\limits_{t\in [0,T/u]_{\delta/u}
,s\in I(t_u)}Z_u(t,s)>m(u)}.
}
\end{lem}

In the next lemma we prove that the discrete Pickands constant appearing
in Theorems \ref{TH1}, \ref{TH2} and  \ref{THinf} is well defined, positive and finite.

\begin{lem}\label{lemma_finite_discrete_pick_const} For any $\delta\ge 0$
and $\eta$ a centered Gaussian process with stationary increments, a.s.
continuous sample paths and variance satisfying \textbf{A} it holds, that
\bqn{\label{theo1_case_ii_finit_pick_const}
\limit{S}\frac{\mathcal H_{\eta}([0,S]_\delta)}{S} =
\mathcal H^\delta_{\eta}\in (0,\IF).
}
\end{lem}

The lemma below allows us to give upper bounds for the double-sum terms
appearing in case $\varphi = 0$ in Theorems \ref{TH1} and \ref{THinf}.
\begin{lem}\label{lemma_case_i)} Assume that $\varphi =0$. Then
uniformly for $t\neq s\in I(t_u)$ and all
large $u$ with some $\ve>0$ it holds that
\bqny{
\pk{X_u(t)>m(u),X_u(s)>m(u)} \le u^{-1/2-\ve}\Psi(m(u)).
}
\end{lem}

The proofs of Lemmas \ref{the_main_area}, \ref{lemma_finite_discrete_pick_const}
and \ref{lemma_case_i)} are given in the Appendix.
\\

\subsection{\prooftheo{TH1}.} Taking $T=0$ in Lemma
\ref{the_main_area} we obtain that
\bqn{\label{theo1_mainint}
\pk{Q_{\delta_X}(0)>u} \sim \pk{\sup\limits_{t\in I(t_u)}X_u(t)>m(u)}, \ \ \ u \to
\IF
,}
where $I(t_u)$ is defined in \eqref{I_t_u_def}.
Next we consider
3 cases:  $\varphi = 0$,  $\varphi \in (0,\IF)$
and $\varphi = \IF$.
\\
\\
{\underline{\emph{Case $\varphi = 0$.}}}
We have by Bonferroni inequality
\bqn{\label{bonf_theo_1.1}\notag
 \sum\limits_{t\in I(t_u)}\pk{X_u(t)>m(u)} &\ge&
\pk{\sup\limits_{t\in I(t_u)}X_u(t)>m(u)}
\\&\ge&
\sum\limits_{t\in I(t_u)}\pk{X_u(t)>m(u)}-
\sum\limits_{t\neq s \in I(t_u)}\pk{X_u(t)>m(u),X_u(s)>m(u)}.
}
There are less then $\mQ u\ln^2 u$ summands in the double-sum above, hence
by Lemma \ref{lemma_case_i)}
 we have
\bqn{\label{double_sum}
\sum\limits_{t\neq s \in I(t_u)}\pk{X_u(t)>m(u),X_u(s)>m(u)}
\le \mQ \ln^2(u) u^{1/2-\ve'}\Psi(m(u)), \ \ \ \ u \to \IF.
}
Next we focus on calculation of the single sum in \eqref{bonf_theo_1.1}.
Since by Lemma \ref{L1} $\sup\limits_{t\in I(t_u)}|\sigma_{X_u}
(t)-1| \to 0$ as $u \to \IF$
the following inequality (see, e.g., Lemma 2.1 in \cite{PicandsA})
$$(1-\frac{1}{x^2})\frac{1}{\sqrt{2\pi}x}e^{-x^2/2}
\le \Psi(x)
\le \frac{1}{\sqrt{2\pi}x}e^{-x^2/2}, \ \ \ \ x>0,
$$
implies as $ u \to \IF$
\bqn{\label{proof_alpha<1/2_sing_sum}
\sum\limits_{t\in I(t_u)}\pk{X_u(t)>m(u)}
\notag&=&
\sum\limits_{t \in I(t_u)}\Psi(\frac{m(u)}{\sigma_{X_u}(t)})
\\&\sim&
\sum\limits_{t \in I(t_u)}\frac{\sigma_{X_u}(t)}{\sqrt{2\pi}m(u)}
e^{-\frac{m^2(u)}{2\sigma^2_{X_u}(t)}}
\notag\\&\sim&
\frac{e^{-\frac{m^2(u)}{2}}}{\sqrt{2\pi}m(u)}\!
\sum\limits_{t \in I(t_u)}e^{-\frac{m^2(u)}{2\sigma_{X_u}^2(t)}+
\frac{m^2(u)}{2}}
\notag\\&\sim&
\Psi(m(u))\sum\limits_{t \in I(t_u)}e^{-\frac{m^2(u)}{2}
\frac{(1-\sigma_{X_u}^2(t))}{\sigma^2_{X_u}(t)}}
.}

By Lemma \ref{L1} we have that as $u \to \IF$ the last sum above is
asymptotically equivalent to
\bqn{
\sum\limits_{t \in I(t_u)}e^{-m^2(u)\frac{B}{2A}(t-t_u)^2 }
\notag&=&
\sum\limits_{t \in (-\frac{\ln u}{\sqrt u},\frac{\ln u}
{\sqrt u})_{\delta/u}}e^{-m^2(u)\frac{B}{2A}t^2}
\\ &=&
\frac{u}{\delta m(u)}
\Big(\frac{\delta m(u)}{u}
\sum\limits_{t \in (-\frac{m(u)\ln u }{\sqrt u},\frac{m(u)\ln u }{\sqrt u})
_{\delta m(u)/u}}
e^{-\frac{B}{2A}t^2}\Big)
\notag\\&\sim& \label{sum_converges_to_integral}
\frac{u}{\delta m(u)}\int\limits_\R e^{-\frac{B}{2A}t^2}dt
\\&=&\notag
\frac{u}{\delta m(u)}\sqrt\frac{2\pi A}{B}
=\notag
\frac{u}{\delta m(u)}\frac{\sqrt{2\pi\alpha}}{c(1-\alpha)^{3/2}}
,}
where the asymptotic equivalence in \eqref{sum_converges_to_integral}
holds, since by \eqref{asympt_m(u)},
$\frac{\delta m(u)}{u} \to 0$ and $\frac{m(u)\ln u }{\sqrt u} \to \IF$ as
$u \to \IF$.
Thus,
\bqn{\label{theo1_case_i_finalsum}
\sum\limits_{t\in I(t_u)}\pk{X_u(t)>m(u)}  \sim
\frac{\sqrt{2\pi\alpha}u\Psi(m(u))}{\delta c(1-\alpha)^{3/2} m(u)},
\ \ \ u \to \IF
}
and hence by \eqref{asympt_m(u)}, \eqref{bonf_theo_1.1} and
\eqref{double_sum} we have that
$$\pk{\sup\limits_{t\in I(t_u)}X_u(t)>m(u)} \sim
\frac{\sqrt{2\pi\alpha}u\Psi(m(u))}{\delta c(1-\alpha)^{3/2} m(u)},
\ \ \ u \to \IF
$$
and the claim follows by \eqref{theo1_mainint}.
\\

{\underline{\emph{Cases $\varphi\in (0,\infty)$ and $\varphi=\infty$.}}}
For any fixed $u>0$ and $S \in \{0,\delta,2\delta,...\}$ denote
$$N_u = \lceil \frac{ u \delta_u}{S\Delta(u)}\rceil, \ \ \ \
t_j = \frac{\Delta(u)jS}{u}, \ \ \ \
\Delta_{j,S,u} = [t_u+t_j,t_u+t_{j+1}]_{\delta/u},\ \ \ \
j\in [-N_u-1,N_u],$$
where $\lceil\cdot\rceil$ is the ceiling function.
We have by Bonferroni inequality that
\bqn{\label{theo1_case_ii_doublesum}
 \sum\limits_{-N_u\le j \le N_u-1}p_{j,S,u}-
\sum\limits_{-N_u-1\le i\neq j \le N_u}p_{i,j;S,u}
\le
\pk{\sup\limits_{t\in I(t_u)}X_u(t)>m(u)}
\le \sum\limits_{-N_u-1\le j \le N_u}p_{j,S,u},
}
where
\bqny{
p_{j,S,u} = \pk{\sup\limits_{t\in \Delta_{j,S,u}} X_u(t)>m(u) }
\ \ \ \ \text{ and } \ \ \ \
p_{i,j;S,u} = \pk{\sup\limits_{t \in \Delta_{j,S,u}} X_u(t)>m(u),
\sup\limits_{t \in \Delta_{i,S,u}} X_u(t)>m(u)}
.}

By \cite{KrzysPeng2015}
(\emph{proof of lower bound of} $\pi_{T_u}(u)$, p. 288)
we have as $u \to \IF$ and then $S \to \IF$
\[\sum\limits_{-N_u-1\le i\neq j \le N_u}p_{i,j;S,u}
=o\left(\frac{u}{m(u)\Delta(u)}\Psi(m(u))\right).\]
Hence from
the asymtotics of $\sum\limits_{-N_u\le j \le N_u-1}p_{j,S,u}$ given in
\eqref{theo1_case_ii_asympt_doublesum} and
\eqref{theo_1_case_iii_asympt_singlesum} we obtain that
as $u \to \IF$ and then $S \to \IF$
$$\pk{\sup\limits_{t\in I(t_u)}X_u(t)>m(u)}
\sim \sum\limits_{-N_u \le j \le N_u}p_{j,S,u}$$
and we need to calculate the asymptotics of the sum above.
That can be done via uniform approximation of $p_{j,S,u}$ for all
$-N_u-1\le j \le N_u$ separately for cases $\varphi\in (0,\IF)$
and $\varphi=\IF$.
\\

{\underline{\emph{Case $\varphi\in (0,\infty)$.}}}
Let
$\Delta(u) = 1$, $N_u = \lceil \frac{ \sqrt u \ln u}{S}\rceil$, $t_j
= \frac{jS}{u}$ and $\Delta_{j,S,u} = [t_u+t_j,t_u+t_{j+1}]_{\delta/u}$.
We have for any $\ve>0, \ 0\le j\le N_u$ for all $u$
large enough
with $m^-_j(u) = \frac{m(u)}{1-(1-\ve)\frac{B}{2A}(\frac{j S}{u})^2}$
\bqny{
p_{j,S,u} =
\pk{\exists t\in \Delta_{j,S,u}: \overline X_u(t)>
\frac{m(u)}{\sigma_u(t)} }
&\le&
\pk{\exists t\in \Delta_{j,S,u}: \overline X_u(t)>
\frac{m(u)}{1-(1-\ve)\frac{B}{2A}(t-t_u)^2}}
\\&\le&
\pk{\sup\limits_{t\in \Delta_{j,S,u}} \overline X_u(t)>m^-_j(u)}
\\&=&
\pk{\sup\limits_{t\in [0,S]_\delta}
\overline X_u(\frac{t+ut_u+ut_j}{u})>m_j(u)}
\\&=:&
\pk{\sup\limits_{t\in [0,S]_\delta}
\overline X_u'(t)>m_j(u)}
.}
By Lemma \ref{L1} we have that 
$1-r_{\overline X'_u}(t,s) \sim \frac{\sigma^2(|t-s|)}{2\sigma^2(ut_*)},
\ u\to \IF, \ t,s \in [0,S]$.
Thus, by Lemma 1 in \cite{DeK} the last probability above
is asymptotically equal  to
$\mathcal H_{\eta'}([0,S]_\delta) \Psi(m^-_j(u)),$ as $u \to \IF$,
where $\eta'$ is a centered Gaussian process with stationary
increments, a.s. continuous sample paths and variance
(asymptotics of $m(u)$ is given in \eqref{asympt_m(u)})
$$\sigma^2_{\eta'}(t) =
\limit{u}\frac{m^2(u)}{2\sigma^2(ut_*)}\sigma^2(t)
=\limit{u} \frac{u^2(1+ct_*)^2}{2t_*\sigma^2(ut_*)\sigma^2(u)
}\sigma^2(t)
= \frac{2c^2}{\varphi^2}\sigma^2(t).
$$
Note that $\eta'$ and $\eta$ defined in \eqref{eta} have the same distributions.
Thus,
$$
\sum\limits_{0\le j \le N_u}p_{j,S,u}
\le
\mathcal H_{\eta}([0,S]_\delta)\sum\limits_{0\le j \le N_u}
\Psi(m^-_j(u)), \ \ \ u \to \IF.
$$
Next, as $u \to \IF$, (set $C_- = \frac{(1-\ve)B}{2A}$)
using the same reasoning as in case $\varphi=0$, we have
\bqny{
\frac{\sum\limits_{0\le j \le N_u}
\Psi(m^-_j(u))}{\Psi(m(u))}
&\sim&
\sum\limits_{0\le j \le N_u}
e^{-\frac{m^2(u)}{2}( \frac{1}{(1-C_-(\frac{j S}{u})^2)^2}
-1)}
\\&\sim&
\sum\limits_{0\le j \le N_u}
e^{-\frac{m^2(u)}{2}2C_-(\frac{j S}{u})^2}
\\&=&
\sum\limits_{ \frac{S m(u)}{u}
\in [0,\frac{N_u S m(u)}{u}]_{\frac{S m(u)}{u}}}
e^{-C_-(\frac{jSm(u)}{u})^2}
\\&=&
\frac{u}{S m(u)}\Big(\frac{S m(u)}{u}
\sum\limits_{t \in [0,\frac{m(u)\ln u }{\sqrt u}]_{\frac{S m(u)}{u}}}
e^{-C_-t^2}\Big).
}
Due to \eqref{asympt_m(u)} we have $\frac{u}{S m(u)} \to 0$ and
$\frac{m(u)\ln u }{\sqrt u} \to \IF$, as $u \to \IF$. Thus, the sum
above converges to
$\int\limits_0^\IF e^{-C_-t^2}dt=\frac{\sqrt\pi}{2\sqrt C_-}$
as $u \to \IF$.
Similar calculation can be done for $j<0$. Hence, 
we have
as $u \to \IF$ and then $ S \to \IF$
$$\sum\limits_{-N_u\le j \le N_u}p_{j,S,u} \le
\mathcal H_{\eta}([0,S]_\delta)
\frac{u\Psi(m(u))}{S m(u)}\frac{\sqrt\pi}{\sqrt C_-}. $$

By Lemma \ref{lemma_finite_discrete_pick_const} we know that
$\frac{\mathcal H_{\eta}([0,S]_\delta)}{S} \to
\mathcal H^\delta_{\eta}\in (0,\IF)$  as $S \to \IF$. Hence,
letting $S \to \IF$, we have
$$\sum\limits_{-N_u\le j \le N_u}p_{j,S,u} \le
\mathcal H^\delta_{\eta}
\frac{u\Psi(m(u))}{ m(u)}\frac{\sqrt\pi}{\sqrt C_-}(1+o(1)), \ \ \ \
u \to \IF. $$
By the same arguments we have the asymptotic lower bound
$$\sum\limits_{-N_u\le j \le N_u}p_{j,S,u} \ge
\mathcal H^\delta_{\eta}
\frac{u\Psi(m(u))}{ m(u)}\frac{\sqrt\pi}{\sqrt C_+}(1+o(1)), \ \ \ \
u \to \IF,$$
with $C_+ = \frac{(1+\ve)B}{2A}$.
Hence letting $\ve \to 0$ we have that, as $S \to \IF$ and then $u \to \IF$,
\bqn{\label{theo1_case_ii_asympt_doublesum}
\sum\limits_{-N_u\le j \le N_u}p_{j,S,u} \sim
\mathcal H^\delta_\eta \frac{u}{m(u)}
\sqrt{\frac{2\pi A}{ B}}\Psi(m(u)),}
which completes the proof of the
case $\varphi\in(0,\infty)$.
\\

\underline{Case $\varphi=\infty$.}
In this case we have that $\Delta(u),N_u \to \IF$ and $\Delta(u) / u \to 0$ as $u \to \IF$.
Moreover, for $m_k(u) = \frac{m(u)}{1-C_+|t_k|^2}, \ \
-N_u\le k \le N_u$ for large $S,u$
\bqny{
p_{k,S,u}
&\ge&
\pk{\sup\limits_{t\in [t_k,t_{k+1}]_{\delta/u}}\overline X_u(t)>m_k(u)}
\\&\ge&
\pk{\sup\limits_{t\in [t_k,t_{k+1}]_{
\frac{\delta[\ve_1\Delta(u)]} {u}
}}
\overline X_u(t)>m_k(u)}
\\&\ge& \pk{\sup\limits_{t\in [0,S-\frac{\delta[\ve_1\Delta(u)]}
{\Delta(u)}]_{\delta\ve_1\frac{[\ve_1\Delta(u)]}
{\ve_1\Delta(u)}}} \overline X_u(\frac{t\Delta(u)}{u}+t_k)>m_k(u)}
\\&\ge& \pk{\sup\limits_{t\in [0,S(1-\ve_2)]
_{\delta\ve_1}} \overline X_u
(\frac{t[\Delta(u)\ve_1]}{u\ve_1}+t_k)>m_k(u)}
\\&=:& \pk{\sup\limits_{t\in [0,S(1-\ve_2)]
_{\delta\ve_1}} Y_u(t)>m_k(u)}
, \ \ \ u \to \IF,
}
where $\ve_1,\ve_2$ are any small positive numbers.
As follows from Lemma \ref{L1}, correlation of $Y_u$,
as $u \to \IF$, has expansion
\bqny{r_{Y_u}(t,s) = 1-\frac{\sigma^2(\Delta(u))|t-s|^{2\alpha}}
{2\sigma^2(u)t_*^{2\alpha}}+
o\left(\frac{\sigma^2(\Delta(u))|t-s|^{2\alpha}}
{\sigma^2(u)}\right)
.}
Next by Lemma 5.1 in \cite{KrzysPeng2015} using the notation introduced therein with
index set $K$ consisting of 1 element and
$$g(u) = m_k(u),\quad \theta(u,s,t) = |t-s|^{2\alpha},
\quad V = B_{\alpha}, \quad \sigma_V(t)=|t|^{2\alpha}$$
we have uniformly for $-N_u\le k\le N_u$
\bqny{
\pk{\sup\limits_{t\in [0,(1-\ve_2)S]_{\delta\ve_1}}Y_u(t)>m_k(u)} \sim
\mathcal{H}_{B_{\alpha}}([0,S(1-\ve_2)]_{\delta\ve_1})\Psi(m_k(u)), \ \ \ u \to \IF.
}
Thus, for large $S$, as $u \to \IF$,
\bqny{
\sum\limits_{-N_u\le k\le N_u}p_{k,S,u}
\ge \mathcal{H}_{B_{\alpha}}([0,(1-\ve_2)S]_{\delta\ve_1})
\sum\limits_{-N_u\le k\le N_u}\Psi(m_k(u))(1+o(1))
.}
Next we calculate the sum above. Similarly to the previous cases we have
as $u \to \IF$ with
$\widehat{C}_+ = \frac{C_+(1+ct_*)^2}{(t_*)^{2\alpha}}$ and $l_u =
\frac{S\Delta(u)}{\sigma(u)} \to 0$, that
\bqny{
\frac{\sum\limits_{k=- N_u}^{N_u}\Psi(m_k(u))}{\Psi(m(u))}
&\sim&
\sum\limits_{k=- N_u}^{N_u}e^{-\frac{m^2(u)}{2}(\frac{1}{(1-
C_+t_k^2)^2}-1)}
\\&\sim&
\sum\limits_{ k=- N_u}^{N_u}
e^{-C_+m^2(u)t_k^2}
\sim
\sum\limits_{k=-\frac{ \sigma(u)\ln u}{S\Delta(u)}}^
{\frac{ \sigma(u)\ln u}{S\Delta(u)}}
e^{-\frac{C_+(1+ct_*)^2}{t_*^{2\alpha}}
(\frac{kS\Delta(u)}{\sigma(u)})^2}
\\&=&
\sum\limits_{kl_u\in(-\ln u,\ln u)_{l_u}}
e^{-\widehat{C}_+(kl_u)^2}
=
\frac{1}{l_u}\Big(l_u
\sum\limits_{t\in(-\ln u,\ln u)_{l_u}}
e^{-\widehat{C}_+t^2}\Big).
}
Since $l_u \to 0$ as $u \to \IF$, the expression in the parentheses above
converges as $u \to \IF$ to
\bqny{
\int\limits_\R e^{-\widehat{C}_+t^2}dt
=\frac{\sqrt \pi}{\sqrt{\widehat{C}_+}}
=\sqrt{\frac{2A\pi}{B}}
\frac{t_*^{\alpha}}{1+ct_*}\frac{1}{\sqrt{1+\ve}}
.}
Thus, summarizing the calculations above
we have as $u \to \IF$ for large $S$
\bqny{
\sum\limits_{-N_u\le k\le N_u}
p_{k,S,u}
&\ge&
\frac{1}{\sqrt{1+\ve}}\frac{1}{S}\mathcal{H}_{B_{\alpha}}([0,
(1-\ve_2)S]_{\delta\ve_1})\Psi(m(u))
\frac{\sigma(u)}{\Delta(u)}\sqrt{\frac{2A\pi}{B}}
\frac{t_*^{\alpha}}{1+ct_*}(1+o(1))
.}
Letting $S \to \IF$ then $\ve_2 \to 0$ and then $\ve_1\to 0$, in view of
Lemma 12.2.7 ii) and Remark 12.2.10 in \cite{LeadbetterExtremes}
we obtain
$\frac{1}{S}\mathcal{H}_{B_{\alpha}}([0,(1-\ve_2)S]_{\delta\ve_1})
\to \mathcal H_{B_{\alpha}}$.
Letting then $\ve \to 0$ we get the asymptotic lower bound
\bqn{\label{theo_1_case_iii_asympt_singlesum}
\sum\limits_{-N_u\le k\le N_u}
p_{k,S,u}
&\ge&
\mathcal H_{B_{\alpha}} \Psi(m(u))
\frac{\sigma(u)}{\Delta(u)}\sqrt{\frac{2A\pi}{B}}
\frac{t_*^{\alpha}}{1+ct_*}(1+o(1)), \ \ \ u \to \IF.
}
Similarly,
we have 
for $\ve>0$
with $m_k^-(u) = \frac{m(u)}{1-C_-|t_k|^2}$ (recall,
$C_- = \frac{B(1-\ve)}{2A}$) as
$S\to \IF$
\bqny{
\sum\limits_{-N_u\le k\le N_u}
p_{k,S,u}
&\le&
\frac{1}{\sqrt{1-\ve}}\mathcal{H}_{B_{\alpha}}\Psi(m(u))
\frac{\sigma(u)}{\Delta(u)}\sqrt{\frac{2A\pi}{B}}
\frac{t_*^{\alpha}}{1+ct_*}(1+o(1)), \ \ \ \ u \to \IF
.}
Letting $\ve \to 0$ we obtain the right side in
\eqref{theo_1_case_iii_asympt_singlesum}, thus, the claim is established.
\QED
\\

\subsection{\textbf{Proof of Theorem \ref{TH2}.}}
By Lemma \ref{the_main_area} we have
\bqn{\label{theo2_main_area}
\psi_{T,\delta}^{\sup}(u) \sim \pk{\sup\limits_{t\in [0,T/u]_{\delta/u}
,s-t\in I(t_u)}Z_u(t,s)>m(u)}, \ \ \ u \to \IF.
}
As in the proof of Theorem \ref{TH1} next we consider three cases:
$\varphi = 0$, $\varphi \in (0,\IF)$ and $\varphi = \IF$.
\\
\\
{\underline{\emph{Case $\varphi = 0$.}}}
 We have by Bonferroni inequality
\bqny{
\sum\limits_{t\in [0,\frac{T}{u}]_{\frac{\delta}{u}}}
\!\sum\limits_{s\in I(t_u)} \!\!\!\pk{Z_u(t,s)\!>\!m(u)}
\!\!\!&\ge&\!\!\!
\pk{\sup\limits_{t\in [0,\frac{T}{u}]_{\frac{\delta}{u}},s\in I(t_u)}
\!\!\!Z_u(t,s)>m(u)}
\\\!\!\!&\ge&\!\!\!\!\!\!\!
\sum\limits_{t\in [0,\frac{T}{u}]_{\frac{\delta}{u}}}
\!\sum\limits_{s\in I(t_u)}
\!\!\!\pk{Z_u(t,s)\!>\!m(u)}-\!\!\!
\sum\limits_{t\in [0,\frac{T}{u}]_{\frac{\delta}{u}}}\!
\sum\limits_{\overset{s_1, s_2 \in I(t_u)}{s_1\neq s_2}}\!\!\!\!\!
\pk{Z_u(t,s_1),Z_u(t,s_2)\!>\!m(u)}\!.
}

For the last term above we have
\bqn{\label{qqqq}\notag
\sum\limits_{t\in [0,\frac{T}{u}]_{\frac{\delta}{u}}}\!
\sum\limits_{\overset{s_1, s_2 \in I(t_u)}{s_1\neq s_2}}\!\!\!\!\!
\pk{Z_u(t,s_1),Z_u(t,s_2)\!>\!m(u)}
\!\!&\le&\!\!
 (1+[\frac{T}{\delta}])
\sup\limits_{t\in [0,\frac{T}{u}]_{\frac{\delta}{u}}}
\sum\limits_{\overset{s_1, s_2 \in I(t_u)}{s_1\neq s_2}}\!\!\!
\pk{X_u(s_1-t),X_u(s_2-t)>m(u)}
\\\!\!&\le&\!\!
\mQ (1+[\frac{T}{\delta}]) \ln^2(u) u^{1/2-\ve'}\Psi(m(u)),
}
where the last inequality above follows from
\eqref{double_sum}.
Next we have
 \bqny{
\sum\limits_{t\in [0,\frac{T}{u}]_{\frac{\delta}{u}}}
\sum\limits_{s\in I(t_u)} \pk{Z_u(t,s)>m(u)}
&=&
\sum\limits_{t\in [0,\frac{T}{u}]_{\frac{\delta}{u}}}
\sum\limits_{s\in I(t_u)} \pk{X_u(s-t)>m(u)}
 \\&\sim&
(1+[\frac{T}{\delta}])\sum\limits_{\tau \in I(t_u)}
\pk{X_u(\tau)>m(u)}, \ \ \ u \to \IF
.}
The asymptotics of the last sum above was calculated in
\eqref{proof_alpha<1/2_sing_sum},
hence
\bqny{\sum\limits_{t\in [0,\frac{T}{u}]_{\frac{\delta}{u}}}
\sum\limits_{s_1\neq s_2 \in I(t_u)} \pk{Z_u(t,s)>m(u)}
\sim(1+[\frac{T}{\delta}])
\frac{\sqrt{2\pi\alpha}u\Psi(m(u))}{\delta c(1-\alpha)^{3/2} m(u)}, \ \ \  \
u \to \IF.
}
By the line above combined with \eqref{qqqq} we obtain
$$
\pk{\sup\limits_{t\in [0,\frac{T}{u}]_{\frac{\delta}{u}},s\in I(t_u)}
Z_u(t,s)>m(u)}\sim(1+[\frac{T}{\delta}])
\frac{\sqrt{2\pi\alpha}u\Psi(m(u))}{\delta c(1-\alpha)^{3/2} m(u)}, \ \ \  \
u \to \IF
$$
and the claim follows by \eqref{theo2_main_area}.
\\
\\
{\underline{\emph{Case $\varphi \in (0,\IF)$.}}}
With the notation of Theorem \ref{TH1}
we have by Bonferroni inequality for $u>0$
\bqny{
\sum\limits_{j=-N_u}^{N_u}q_{j,S,u}
\ge
\pk{\sup\limits_{t\in [0,\frac{T}{u}]_{\frac{\delta}{u}},s\in I(t_u)} Z_u(t,s)
>m(u)}
\ge\sum\limits_{j=-N_u}^{N_u}q_{j,S,u}
-
\sum\limits_{-N_u \le i<j\le N_u}q_{i,j;S,u},
}
where
\bqny{
q_{j,S,u} = \pk{\!
\sup\limits_{t\in [0,\frac{T}{u}]_{\frac{\delta}{u}},
s \in \Delta_{j,S,u}}\!\!\!\!\!\!\!\!\!Z_u(t,s)\!>\!m(u)\!}\!
, \
q_{i,j;S,u} =
\pk{\!\!
\exists t\!\in\! [0,\!\frac{T}{u}]_{\frac{\delta}{u}}\!\!:\!\!\!
\sup\limits_{s \in \Delta_{j,S,u}}\!\!\!Z_u(t,s)\!>\!m(u),\!\!
\sup\limits_{s \in \Delta_{i,S,u}}\!\!\!Z_u(t,s)\!>\!m(u)\!\!}
\!.}

By \cite{KrzysPeng2015}
(\emph{proof of lower bound of} $\pi_{T_u}(u)$, p. 288)
we have that as $u \to \IF$ and then $S\to \IF$
\bqny{
\sum\limits_{-N_u \le i<j\le N_u}q_{i,j;S,u} =
o\left(\frac{u}{m(u)\Delta(u)}\Psi(m(u))\right).
}
Since
$$\pk{\sup\limits_{t\in [0,\frac{T}{u}]_{\frac{\delta}{u}},s\in I(t_u)} Z_u(t,s)
>m(u)} \ge\pk{\sup\limits_{t\in I(t_u)}X_u(t)>m(u)} \ge
C\frac{u}{m(u)\Delta(u)}\Psi(m(u)),
\quad u \to \IF,
$$
we have that as $u \to \IF$ and then $ S \to \IF$
$$
\pk{\sup\limits_{t\in [0,\frac{T}{u}]_{\frac{\delta}{u}},s\in I(t_u)} Z_u(t,s)
>m(u)} \sim
\sum\limits_{j=-N_u}^{N_u}q_{j,S,u}
$$
and we need to calculate the asymptotics of the sum above.
Next we uniformly approximate each summand in the sum above.
For $\ve>0, j\ge 1, S>T$ and $u$ large enough
we have (recall, $m_{j-1}^-(u) =
\frac{m(u)}{1-(1-\ve)\frac{B}{2A}(\frac{(j-1)S}{u})^2}$)
\bqny{
q_{j,S,u}
&=&
\pk{\!\exists (t,s) \!\in\! [0,\frac{T}{u}]_{\delta/u}\!\times\!
 \Delta_{j,S,u}\!:\!\overline Z_u(t,s)>\frac{m(u)}{\sigma_u(s-t)}}
\\&\le&
\pk{\!\exists (t,s) \!\in\! [0,\frac{T}{u}]_{\delta/u}\!\times\!
 \Delta_{j,S,u}\!:\! \overline Z_u(t,s)\!>\! \frac{m(u)}{1\!-\!
 (1-\ve)\frac{B}{2A}(s\!-\! t\!-\!t_u)^2}}
\\&=&
\pk{\sup\limits_{t\in [0,T]_{\delta},
s \in [0,S]_\delta}\overline Z_u(t/u,t_u+s/u)>\frac{m(u)}{1-C_-
(\frac{(j-1)S}{u})^2}}
\\&=:&
\pk{\sup\limits_{t\in [0,T]_{\delta},
s \in [0,S]_\delta}\overline Z'_u(t,s)>m^-_{j-1}(u)}
.}
By \eqref{r_Z_u_formula}
correlation of $\overline Z_u'$ satisfies
\bqn{\label{correlation_Z'_u}
1-r_{\overline Z'_u}(t,s,t_1,s_1) \sim \frac{\sigma(|s-s_1|)^2+\sigma^2(|t-t_1|)}
{2\sigma^2(ut_*)}, \ \ \ (t,s) \in [0,T]\times [0,S], \ \ \  u \to \IF .
}
Applying Lemma 5.1 in \cite{KrzysPeng2015} with
$$\Phi := \sup, \quad \theta(u,s,t) := \frac{2c^2}{\varphi^2}(\sigma^2(|s-s_1|)
+\sigma^2(|t-t_1|)), \ \ \ \ V(t,s) := \frac{\sqrt 2 c}{\varphi}
(X^{(1)}(t)+X^{(2)}(s)), $$
where $X^{(1)}$ and $X^{(2)}$ are independent copies of $X$ we have
$$
\pk{\sup\limits_{t\in [0,T]_{\delta},
s \in [0,S]_\delta}\overline Z'_u(t,s)>m_{j-1}^-(u)} \sim
\mathcal{H}_{\frac{\sqrt 2 c}{\varphi} X}([0,T]_\delta)
\mathcal{H}_{\frac{\sqrt 2 c}{\varphi} X}([0,S]_\delta)
\Psi(m_{j-1}^-(u))
.$$

Finally, for $\ve>0,j\ge 1,S>T$ and $u$ large we have
$$\pk{\sup\limits_{t\in [0,\frac{T}{u}]_{\delta/u},
s \in \Delta_{j,S,u}}Z_u(t,s)>m(u)} \le
\mathcal{H}_{\frac{\sqrt 2 c}{\varphi} X}([0,T]_\delta)
\mathcal{H}_{\frac{\sqrt 2 c}{\varphi} X}([0,S]_\delta)
\Psi(m_{j-1}^-(u))(1+o(1)).$$

The rest of the proof is the same as in Theorem \ref{TH1}
case $\varphi\in (0,\IF)$, thus, the claim
is established.
\\
\\
{\underline{\emph{Case $\varphi = \IF$.}}}
 By Theorem \ref{TH1} we have
\bqny{
\psi_{T,\delta}^{\sup}(u)\ge \pk{M_\delta>u} \sim
\mathcal{H}_{B_\alpha} f(u)\Psi(m(u))
,\quad u \to \IF.}
By \eqref{theo2_main_area} we have
\bqny{
\psi_{T,\delta}^{\sup}(u)
\le
\pk{\sup\limits_{t\in [0,T/u],s-t\in (-\delta_u+t_u,t_u+\delta_u
)}Z_u(t,s)>m(u)}(1+o(1)), \quad u \to \IF.
}
From the proof of Theorem 3.1 in
\cite{KrzysPeng2015} it follows that the
last probability above does not exceed
$(1+o(1))\mathcal{H}_{B_\alpha} f(u)\Psi(m(u)), \ u \to \IF$. Combining both bounds
above we obtain the claim. \QED
\\

\textbf{Proof of Theorem \ref{THinf}.}
{\underline{\emph{Case $\varphi \in (0,\IF)$.}}}
By Lemma \ref{the_main_area} we have
$$\psi_{T,\delta}^{\inf}(u)\sim
\pk{\inf\limits_{t\in [0,\frac{T}{u}]_{\delta/u}}
\sup\limits_{s\in I(t_u)} Z_u(t,s)
>m(u)}, \ \ \ \ u \to \IF.$$

With notation of Theorem \ref{TH1}
repeating the proof of Theorem \ref{TH2} we have as $u \to \IF$
and then $S \to \IF$
$$\pk{\inf\limits_{t\in [0,\frac{T}{u}]_{\delta/u}}
\sup\limits_{s\in I(t_u)} Z_u(t,s)
>m(u)}\sim \sum\limits_{j=-N_u}^{N_u}\pk{
\inf\limits_{t\in [0,\frac{T}{u}]_{\delta/u}}
\sup\limits_{s \in \Delta_{j,S,u}}Z_u(t,s)>m(u)
}.
$$
Next we uniformly approximate each summand in the sum above.
For $\ve>0, j\ge 1, S>T$ and $u$ large enough
similarly to the proof of Theorem \ref{TH2} we obtain
(recall, $m_{j-1}^-(u) =
\frac{m(u)}{1-(1-\ve)\frac{B}{2A}(\frac{(j-1)S}{u})^2}$)
\bqny{
\pk{
\inf\limits_{t\in [0,\frac{T}{u}]_{\delta/u}}
\sup\limits_{s \in \Delta_{j,S,u}}Z_u(t,s)>m(u)
}
\le
\pk{\inf\limits_{t\in [0,T]_{\delta}}
\sup\limits_{s \in [0,S]_\delta}\overline Z'_u(t,s)>m_{j-1}^-(u)}
,}
where $\overline Z_u'(t,s)$ is a centered
Gaussian process with unit variance and correlation
satisfying \eqref{correlation_Z'_u}.
Applying Lemma 5.1 in \cite{KrzysPeng2015} with
(with $X^{(1)},X^{(2)}$ being independent copies of $X$)
$$
\Phi := \inf\sup, \quad
\theta(u,s,t) := \frac{2c^2}{\varphi^2}(\sigma^2(|s-s_1|)
+\sigma^2(|t-t_1|)), \ \ \ \ V(t,s) := \frac{\sqrt 2 c}{\varphi}
(X^{(1)}(t)+X^{(2)}(s)) $$
we have as $u \to \IF$
$$
\pk{\inf\limits_{t\in [0,T]_{\delta}}
\sup\limits_{s \in [0,S]_\delta}\overline Z'_u(t,s)>
m_{j-1}^-(u)} \sim
\E{\inf\limits_{t \in [0,T]_\delta}\sup\limits_{s\in [0,S]_\delta}
e^{\sqrt 2V(t,s)-\var(V(t,s)) } }
\Psi(m_{j-1}^-(u))
.$$
Since $X^{(1)}(t)$ and $X^{(2)}(s)$ are independent we have
$$\E{\inf\limits_{t \in [0,T]_\delta}\sup\limits_{s\in [0,S]_\delta}
e^{\sqrt 2V(t,s)-\var(V(t,s)) } } =
\mathcal{G}_{\frac{\sqrt 2 c}{\varphi} X}([0,T]_\delta)
\mathcal{H}_{\frac{\sqrt 2 c}{\varphi} X}([0,S]_\delta)
.$$
Thus, as $u \to \IF$
\bqny{
\pk{
\inf\limits_{t\in [0,\frac{T}{u}]_{\delta/u}}
\sup\limits_{s \in \Delta_{j,S,u}}Z_u(t,s)>m(u)} \le
\mathcal{H}_{\frac{\sqrt 2 c}{\varphi} X}^{\inf}([0,T]_\delta)
\mathcal{H}_{\frac{\sqrt 2 c}{\varphi} X}([0,S]_\delta)\Psi(m_{j-1}^-(u))
.}
The rest of the proof is the same as in the proof of Theorem \ref{TH1}. Thus, the claim
is established.
\\
\\
{\underline{\emph{Case $\varphi = \IF$.}}}
Let $R(s,t) = X(s)-X(t)-c(s-t), \ t,s \ge 0.$
Using the idea from \cite{DeK}, (the next equation after (2))
we write
\bqny{
\frac{\psi^{\inf}_{T,\delta}(u)}{\psi^{\sup}_{T,\delta}(u)}
&=&
\pk{\inf\limits_{t\in [0,T]_\delta}
\sup\limits_{t\le s\in G_\delta}
R(s,t)>u\Big|\sup\limits_{t\in[0,T]_\delta}
\sup\limits_{t\le s\in G_\delta} R(s,t)>u }
\\&\ge&
1-\sum\limits_{t\in [0,T]_\delta}\left(1-\frac{\pk{
\sup\limits_{t\le s\in G_\delta}R(s,t)>u}}
{\pk{\sup\limits_{a\in [0,T]_\delta}
\sup\limits_{a\le b \in G_\delta }R(a,b)
>u}}\right).
}
By Theorems \ref{TH1} and \ref{TH2} we have that the right part
of the expression above tends to 1 as $u \to \IF$. Thus, we have
$\psi^{\inf}_{T,\delta}(u) \sim \psi_{T,\delta}^{\sup}(u) \sim
\mathcal{H}_{B_\alpha} f(u)\Psi(m(u)), \ u \to \IF$ and the claim follows.
\QED
\\

\textbf{Proof of Proposition \ref{proposition1}.}
Since $T \ge \delta$ with any $K\in [\delta,T]_\delta$,
$a\in (0,t_*), b>0$ and $J(t_u) = [-a+t_u,t_u+b]$
we have
\bqny{\psi^{\inf}_{\delta,T}(u) &\le&
\psi^{\inf}_{\delta,K}(u)
\\&\le&
 \pk{\exists s_1,s_2 \in J(t_u)\cap G_{\frac{\delta}{u}}:
 Z_u(0,s_1)>m(u),Z_u(\frac{K}{u},s_2)>m(u)}
 \\& \  & +
 \pk{\exists s\notin J(t_u): Z_u(0,s)>m(u)}
 \\& =:& p_1(u)+p_2(u)
.}
\emph{Estimation of $p_1(u)$.}
Fix some $s_1,s_2 \in J(t_u) \cap G_{\frac{\delta}{u}}$
and let $(W_1,W_2) =(Z_u(0,s_1),Z_u(\frac{K}{u},s_2)).$
We have that $(W_1,W_2)$ is a centered Gaussian vector with
$Var(W_1),Var(W_2) \le 1$ and correlation
$r_{W_u}(s_1,s_2)$ satisfying
(see \eqref{r_Z_u_formula})
$$1-r_{W_u}(s_1,s_2) \ge \frac{\sigma^2(u|s_1-s_2|)+\sigma^2(K)}
{2\sigma^2(ut_*)}(1+o(1)) \ge (1+o(1))
\frac{\sigma^2(K)}{2t_*^{2\alpha}\sigma^2(u)}, \ \ \ \ \
s_1,s_2 \in J(t_u).
$$
Thus, by Lemma 2.3 in \cite{PicandsA}
\bqny{
\Psi(m(u))^{-1}\pk{W_1>m(u),W_2(u)>m(u)} &\!\!\le\!\!& 3\Psi \big(m(u)\sqrt{
\frac{1-r_{W_u}(s_1,s_2)}{2}}\big)
\\&\!\!=\!\!& 3\Psi\Big((1+o(1))\frac{\sigma(K)m(u)}
{2t_*^\alpha\sigma(u)} \Big)
\!=\! 3\Psi\Big(\!(1+o(1))
\frac{\sigma(K)(1+ct_*)u}{2t_*^{2\alpha}\sigma^2(u)}\Big)
.}

Note that
$\varphi = 0$ implies $\frac{u}{\sigma^2(u)} \to \IF$ as $u \to \IF$.
Thus, since there are less than $Cu^2\ln^2 u$ points in $(J(t_u)\cap
G_{\frac{\delta}{u}}) \times (J(t_u)\cap G_{\frac{\delta}{u}})$
we have with any $ \mQ_K <
\frac{\sigma(K)(1+ct_*)}{2t_*^{2\alpha}}$ by
\eqref{additional_condition_proposition} as $u \to \IF$
\bqny{
p_1(u)&\le& Cu^2(\ln^2 u)\cdot 3\Psi(m(u))\Psi\Big((1+o(1))
\frac{\sigma(K)(1+ct_*)u}{2t_*^{2\alpha}\sigma^2(u)}\Big)
\le \Psi(m(u))\Psi\Big(\mQ_K
\frac{u}{\sigma^2(u)}\Big).
}
\emph{Estimation of $p_2(u)$.}
Since $Z_u(0,s) \overset{d}{=} X_u(s), \ s\ge 0$,
it follows from the estimation of $R_1(u)$ and $R_2(u)$ in
the proof of Lemma \ref{the_main_area} below
(see \eqref{R_1(u)} and \eqref{R_2(u)}, respectively) that
for appropriately chosen $a\in (0,t_*),b>0$ and small
$\ve >0$
$$p_2(u) \le \pk{\exists s\in [0,a]: X_u(s)>m(u)}+
\pk{\exists s\in [t_*+b,\IF): X_u(s)>m(u)} \le
\Psi(m(u))\mQ e^{-\frac{m^2(u)}{u^{2\ve}}}.$$
 Combing this inequality with the upper bound of $p_1(u)$
 we obtain that for any $K\in [\delta,T]_\delta$ it holds that
 $$\psi^{\inf}_{\delta,T}(u) \le \Psi(m(u))
 \Psi\Big(\mQ_K \frac{u}{\sigma^2(u)}\Big), \ \ \ u \to \IF$$
 and taking the supremum with
 respect to $K$ we obtain the claim. \QED
\\
\\
\textbf{Proof of Corollary \ref{korr_integrated_simple}.}
We start with the proof of the first statement. Since
$Z(t)$ is a Gaussian process with stationary increments, a.s. sample paths
and variance satisfying \textbf{A}
applying Theorem \ref{TH1} with parameters
\bqny{
\varphi = \limit{u}\frac{\sigma_\zeta^2(u)}{u} = \frac{2}{G}\in (0,\IF), \
\alpha = 1/2, \ t_* = 1/c, \ \Delta(u) = 1, \ A = 2\sqrt C, \
D = c^2\sqrt c /2
,\\
\ m(u) = \sqrt {2Gc}\sqrt u + \frac{c^{3/2}G^{3/2}D}{\sqrt 2}u^{-1/2}
+o(u^{-1/2}), \ f(u) =  \frac{2\sqrt\pi}{c\sqrt{cG}}\sqrt u+O(u^{-1}), \
\eta(t) = \frac{cG}{\sqrt 2}Z(t)
}
we have $\Psi(m(u))\sim e^{-ucG-c^2G^2D}\frac{1}{2\sqrt{\pi Gcu}}, u \to \IF$
implying
\bqny{
\pk{\exists t\in [0,\IF):Z(t)-ct>u}
 \sim \mathcal H^\delta_{cGZ(t)/\sqrt 2}\frac{1}{c^2G}
 e^{-ucG-c^2G^2D} ,  \ \ \ \ \ u
\to \IF
}
and the first claim follows.
Applying Theorems \ref{TH2} and \ref{THinf} with the same parameters we
obtain the second and third claims, respectively.
\QED

\textbf{Proof of Lemma \ref{the_main_area}.}
We start with the first claim. We have
\bqn{\label{main_inequality}\notag
\pk{\!\sup\limits_{t \in [0,\frac{T}{u}]_{\frac{\delta}{u}}}
\!\sup\limits_{s\in I(t_u)}\!\!
Z_u(t,s) \!>\!m(u) \!\!}\!\!\!&\le&\!\!\!
\psi_{\delta,T}^{\sup}(u)
\\\!\!\!&\le&\!\!\!
\pk{\!\sup\limits_{t \in [0,\frac{T}{u}]_{\frac{\delta}{u}}}
\!\sup\limits_{s\in I(t_u)}\!\!
Z_u(t,s) \!>\!m(u) \!\!}\!+\!
\pk{\!\sup\limits_{t \in [0,\frac{T}{u}]_{\frac{\delta}{u}}}
\!\sup\limits_{s\in (G_{\frac{\delta}{u}}\backslash  I(t_u))}
\!\!Z_u(t,s) \!>\!m(u)\!\! }.
}
Our first aim is to show that
\bqn{\label{negligible_area}
\pk{\sup\limits_{t \in [0,\frac{T}{u}]_{\frac{\delta}{u}}}
\sup\limits_{s\in (G_{\frac{\delta}{u}}\backslash  I(t_u))}
Z_u(t,s) >m(u) }= o(\Psi(m(u))), \ \ \  \  \
u \to \IF.
}

Since for any fixed $0\le t\le s$
random variables $Z_u(t,s)$ and $X_u(s-t)$
have the same distributions we have
with $I'(t_u) = (-\frac{\delta_u}{2}+t_*,\frac{\delta_u}{2}+t_*)\cap
G_{\delta/u}$
\bqny{
\pk{\sup\limits_{t \in [0,T/u]_{\delta/u}}
\sup\limits_{s\in (G_{\delta/u}\backslash  I(t_u))}
Z_u(t,s) >m(u) } &\le&
\sum\limits_{t \in [0,T/u]_{\delta/u}}
\pk{\sup\limits_{s\in (G_{\delta/u}\backslash  I(t_u))}
Z_u(t,s) >m(u)}
\\&=&
\sum\limits_{t \in [0,T/u]_{\delta/u}}
\pk{\sup\limits_{s\in (G_{\delta/u}\backslash  I(t_u))}
X_u(s-t) >m(u)}
\\&\le&
(1+[\frac{T}{\delta}])
\pk{\sup\limits_{s\in (G_{\delta/u}\backslash   I'(t_u))}
X_u(s) >m(u)}.
}

We have that for any chosen small $\ve$ and large $M$ the last
probability above does not exceed
\bqny{
\sum\limits_{t\in (G_{\frac{\delta}{u}}\backslash I'(t_u))}
\!\!\!\!\!\!\!\!\pk{X_u(t)\!>\!m(u)}
\!\!\!&=&\!\!\!
\sum\limits_{t\in (G_{\delta/u}\backslash I'(t_u))}
\!\!\Psi(\frac{m(u)}{\sigma_{X_u}(t)})
\\\!\!\!&\le&\!\!\!
2\Psi(m(u))\!
\Big(\!\!\sum\limits_{t\in [0,\ve]_{\frac{\delta}{u}}}
\!\!\!\!e^{\!-\frac{m^2(u)}{2}(\frac{1}{\sigma_{X_u}^2\!\!(t)}-1)}
\!\!+\!\!\!
\sum\limits_{t\in [M,\IF)_{\frac{\delta}{u}}}
\!\!\!\!\!\!e^{\!\!-\frac{m^2(u)}{2}(\frac{1}{\sigma_{X_u}^2\!\!(t)}-1)}\!
+\!\!\!\!
\sum\limits_{t\in ([\ve,M]_{\frac{\delta}{u}}\backslash I'(t_u))}
\!\!\!\!\!\!\!\!\!\!\!\!\!e^{\!-\frac{m^2(u)}{2}(\frac{1}{\sigma_{X_u}^2\!\!(t)}-1)}
\Big)
\\&=:& 2\Psi(m(u))(R_1(u)+R_2(u)+R_3(u)).
}
Thus, to establish \eqref{negligible_area} we need to prove that
$R_1(u)+R_2(u)+R_3(u) \to 0$ as $u \to \IF$.
By Lemma \ref{L1}, $t_u$ is unique for large $u$ and we have
\bqny{
\sigma_{X_u}(t) =\frac{\sigma(ut)}{u(1+ct)}m(u) =
\frac{\sigma(ut)}{\sigma(ut_u)}
\frac{1+ct_u}{1+ct}.
}
\emph{Estimation of $R_1(u)$.}
We have for all large $u$ and $t\in [0,\ve]_{\delta/u}$
$$\sigma_{X_u}(t) \le \mQ\frac{\sigma(ut)}{\sigma(ut_u)}.$$
i) Assume that $ut\ge \ln u$. Then with $h$ being a slowly
varying function at $\IF$ and $0<\epsilon<\alpha$
 by Potter's theorem (Theorem 1.5.6 in \cite{BI1989}) we have
$$\frac{\sigma(ut)}{\sigma(ut_u)} = (\frac{t}{t_u})^{\alpha}
\frac{h(ut)}{h(ut_u)} \le \mQ t^{\alpha} (\frac{t_u}{t})^{\epsilon}
\le \mQ t^{\alpha-\epsilon}.$$
ii) Assume that $ut<\ln u$. Since $t\in [0,\ve]_{\delta/u}$ we
have $ut\ge \delta$
for $t\neq 0.$ Then for $\epsilon\in(0,\alpha)$ and large $u$
$$\frac{\sigma(ut)}{\sigma(ut_u)} \le
u^{-(\alpha-\epsilon)}
\sup\limits_{t\in [\delta,\ln u]}\sigma(t)\le
u^{-(\alpha-\epsilon)}\ln u.$$
Combining the above inequalities we have that for sufficiently small
$\ve$ and for all
$t\in [0,\ve]_{\delta/u}$ uniformly for large $u$ it holds that
$\frac{1}{\sigma_{X_u}^2(t)}-1\ge 2$.
Thus, for small enough $\ve>0$
\bqn{\label{R_1(u)}
R_1(u) \le \mQ u e^{-m^2(u)} \to 0, \ \ \ u \to \IF
.}
\emph{Estimation of $R_2(u)$.}
By Potter's theorem we have for $M$ large enough and
$0<\epsilon'<1-\alpha$
\bqny{
\frac{\sigma(ut)}{\sigma(ut_u)} \le \mQ
\left(\frac{t}{t_u}\right)^{\alpha+\epsilon'}
.}
Since $t_u \to t^*$ as $u \to \IF$ we have for some small $\epsilon>0$
$$\sigma_{X_u}(t) \le \mQ \frac{t^{\alpha+\epsilon'}}{1+ct}
\le t^{-\epsilon},$$
hence for all $t>M$ uniformly for $u$ large it holds that
$\frac{1}{\sigma^2_{X_u}(t)}-1 \ge 2t^{\epsilon}.$
Choosing $M$ large enough and $\ve$ sufficiently small
we have as $u \to \IF$
\bqn{\label{R_2(u)}
R_2(u) \le
\sum\limits_{t\in [M,\IF)_{\delta/u}}e^{-t^{\epsilon}m^2(u)}
= \sum\limits_{t\in [M,\IF)_\delta}
e^{-t^{\epsilon}\frac{m^2(u)}{u^\ve}}
\le\mQ e^{-\frac{m^2(u)}{u^\ve}}
  \to 0.
}
\emph{Estimation of $R_3(u)$.}
We have by Lemma \ref{L1} that with some $\mQ>0$
$$1-\sigma_{X_u}(t)\ge \mQ(t-t_u)^2 , \ \ \ t \in [\ve, M]$$
and hence by \eqref{asympt_m(u)}
for $t\in [\ve,M]_{\frac{\delta}{u}}\backslash I(t_u)$ it holds that
$$m^2(u)(\frac{1}{\sigma_{X_u}^2(t)}-1) \ge
m^2(u)(1-\sigma_{X_u}(t))\ge \mQ(t-t_u)^2m^2(u)\ge
\mQ \ln^2 u.$$

Thus, for $t\in [\ve,M]_{\frac{\delta}{u}}\backslash I(t_u)$ it holds that
$e^{-\frac{m^2(u)}{2}
(\frac{1}{\sigma_u^2(t)}-1)} \le u^{-\mQ \ln u},$
and we obtain
\bqny{
R_3(u) \le u\mQ_1  u^{-\mQ_2\ln u}
\to 0, \ \ \ u \to \IF.
}
Combining the estimate above with
\eqref{R_1(u)} and \eqref{R_2(u)}
obtain \eqref{negligible_area}.
It follows from the calculations in Theorem \ref{TH2} that
$$\pk{\sup\limits_{t \in [0,\frac{T}{u}]_{\frac{\delta}{u}}}
\sup\limits_{s\in I(t_u)}
Z_u(t,s) >m(u) }
 \ge \Psi(m(u)), \ \ \ \ u \to \IF$$
and the first claim follows by
\eqref{main_inequality} and \eqref{negligible_area}.
\\\\
Next we show the
second statement of the lemma. Again, by Bonferroni inequality we have
\bqny{
\pk{\!\inf\limits_{t \in [0,\frac{T}{u}]_{\frac{\delta}{u}}}
\!\sup\limits_{s\in I(t_u)}\!\!
Z_u(t,s) \!>\!m(u) \!\!}\!\!\!&\le&\!\!\!
\psi_{\delta,T}^{\inf}(u)
\\\!\!\!&\le&\!\!\!
\pk{\!\inf\limits_{t \in [0,\frac{T}{u}]_{\frac{\delta}{u}}}
\!\sup\limits_{s\in I(t_u)}\!\!
Z_u(t,s) \!>\!m(u) \!\!}\!+\!
\pk{\!\inf\limits_{t \in [0,\frac{T}{u}]_{\frac{\delta}{u}}}
\!\sup\limits_{s\in (G_{\frac{\delta}{u}}\backslash  I(t_u))}
\!\!Z_u(t,s) \!>\!m(u)\!\! }.
}
It follows from the calculations in Theorem \ref{TH2} that
$$\pk{\inf\limits_{t \in [0,\frac{T}{u}]_{\frac{\delta}{u}}}
\sup\limits_{s\in I(t_u)}
Z_u(t,s) >m(u) }
 \ge \Psi(m(u)), \ \ \ \ u \to \IF$$
and by \eqref{negligible_area} we obtain that
\bqny{
\pk{\inf\limits_{t \in [0,\frac{T}{u}]_{\frac{\delta}{u}}}
\sup\limits_{s\in (G_{\frac{\delta}{u}}\backslash  I(t_u))}
Z_u(t,s) >m(u) }= o(\Psi(m(u))), \ \ \  \  \
u \to \IF.
}
Combining both statements above we obtain the second claim of the lemma.
\QED \\

\textbf{Proof of Lemma \ref{lemma_case_i)}.}
Fix some $t\neq s \in I(t_u)$.
Since $\sigma_{X_u}(t),\sigma_{X_u}(s)\le 1$ we have by Lemma 2.3 in
\cite{PicandsA}, with $r_u(t,s) = \corr(\overline X_u(t),
\overline X_u(s))$,
\bqny{
\pk{X_u(t)>m(u),X_u(s)>m(u)}
\le
\pk{\overline X_u(t)>m(u),\overline X_u(s)>m(u)}
\le
\Psi(m(u))\Psi\big(m(u)\sqrt{\frac{1-r_u(t,s)}{2}}\big)
.}
If $\alpha<1/2$, then by Lemma \ref{L1} as $u \to \IF$
for some $\epsilon>0$ it holds that
$m(u)\sqrt{\frac{1-r_u(t,s)}{2}}\ge u^{\epsilon}$
hence uniformly for $t\neq s \in I(t_u)$
\bqn{\label{appendix_theo1_case_i_doublesum_exp_negligib}
\Psi(m(u)\sqrt{\frac{1-r_u(t,s)}{2}})\le e^{-\frac{1}{3}u^{2
\epsilon}}, \ \ \ u \to \IF
.}
If $\alpha = 1/2$, then $t_*= 1/c$ and by Lemma \ref{L1} and
\eqref{theo1_case_i_assumption_sigma}
we have that for some $\ve'>0$ as $u \to \IF$
\bqny{
\frac{m^2(u)}{2}\frac{1-r_u(t,s)}{2} \sim
\frac{u^2(1+ct_*)^2\sigma^2(u|s-t|)}{8\sigma^4(ut_*)}
\sim
\frac{u^2c^2\sigma^2(u|s-t|)}{2\sigma^4(u)}
\ge
(1/2+\ve') \ln u
}
implying as $u \to \IF$
\bqny{
\Psi(m(u)\sqrt{\frac{1-r_u(t,s)}{2}})\le u^{-1/2-\ve'}
.}
Combining the line above with
\eqref{appendix_theo1_case_i_doublesum_exp_negligib} we obtain the claim.
\QED \\\\
\textbf{Proof of Lemma \ref{lemma_finite_discrete_pick_const}.}
Let $a=K\delta$, where $K$ is a large natural number that we
shall choose later on. By the proofs of Theorem 15 and Lemma 16
in \cite{DebickiMichnaRolski_Pickconst} we have with $\sigma^2_\eta$
being the variance of $\eta$ that
$$\liminf\limits_{S\to \IF}\frac{\mathcal H_\eta([0,S]_a)}{S} \ge \frac{1}{a}\Big(1-\frac{2}{a}
\int\limits_0^\IF e^{-\frac{\sigma_\eta^2(t)}{4}}dt \Big).$$
We have that for all $u$ large enough
$\sigma_\eta^2(t) \ge \mQ t $
implying that $\int\limits_0^\IF e^{-\frac{\sigma_\eta^2(t)}{4}}
dt < \IF$. Choosing sufficiently large $K$ we have
$$\liminf\limits_{S\to \IF}\frac{\mathcal H_\eta([0,S]_\delta)}{S} \ge
\liminf\limits_{S\to \IF}\frac{\mathcal H_\eta([0,S]_a)}{S}\ge \frac{1}{a}\cdot \frac{1}{2} >0.$$

Next we prove that for
$I(S):=\frac{\mathcal H_\eta([0,S]_\delta)}{S}$
it holds that for large $S \in G_\delta$
\bqn{\label{Appendix_nonincreasing_pick_const}
I(S) \ge I(S+\delta).
}
We have
\bqny{
(S+\delta)I(S+\delta) \le \E{\sup\limits_{t\in \{0,\delta,...,S\}}
e^{\sqrt 2 \eta(t)-\sigma^2_\eta(t)}}
+\E{e^{\sqrt 2 \eta(S+\delta)-\sigma^2_\eta(S+\delta)}}
F(S+\delta)
=SI(S)+
F(S+\delta),
}
\!\!where
$F(M) = \pk{\underset{t \in \{0,\delta,...,M\}}{\argmax}
(\sqrt 2 \eta(t)-\sigma^2_\eta(t)) = M}$
 for  $M\in G_\delta $.
Thus, to claim \eqref{Appendix_nonincreasing_pick_const}
we need to show that for large $S$
\bqn{\label{equiv}
\delta I(S) \ge F(S+\delta).
}
Since $\liminf\limits_{S\to \IF}I(S)>0$,
we have that $\delta I(S)>\ve$ for all $S$ and some positive $\ve$,
but on the other hand as $S \to \IF$ it holds that
\bqny{
F(S+\delta) \le
\pk{\sqrt 2 \eta(S+\delta)-\sigma^2_\eta(S+\delta) \ge
\sqrt 2 \eta(0)-\sigma^2_\eta(0)}
=
\pk{\sqrt 2 \eta(S+\delta)-\sigma^2_\eta(S+\delta)\ge 0}
\to 0,
}
consequently \eqref{equiv} holds
and hence $I(S)$ is non-increasing for large $S$. Thus,
$\lim\limits_{S \to \IF}I(S) \in (0,\IF)$ and the claim holds. \QED

{\bf Acknowledgement}:
K. D\c{e}bicki
was partially supported by
NCN Grant No 2018/31/B/ST1/00370
(2019-2022).
G. Jasnovidov was supported
by the Ministry of Science and Higher Education of the
Russian Federation, agreement 075-15-2019-1620
date 08/11/2019 and 075-15-2022-289 date 06/04/2022.\\
The Authors would like to thank Professor Enkelejd Hashorva for fruitful
discussions that significantly improved the paper.
\bibliographystyle{ieeetr}
 \bibliography{qu}{}

\end{document}